\documentclass{article}
\usepackage{amsmath, amssymb, amsthm, graphicx,bm,
comment}
\usepackage[all]{xy}

\DeclareMathOperator{\Vol}{Vol}
\DeclareMathOperator{\dVol}{dVol}

\DeclareMathOperator{\grad}{grad} 

\DeclareMathOperator{\area}{area}

\def\e#1\e{\begin{equation}#1\end{equation}}
\def\i#1\i{\begin{itemize}#1\end{itemize}}
\def\ea#1\ea{\begin{align}#1\end{align}}

\def\l{\label}

\theoremstyle{plain}
\newtheorem{thm}{Theorem}[section]
\newtheorem{prop}[thm]{Proposition}
\newtheorem{lem}[thm]{Lemma}

\theoremstyle{definition}

\newtheorem{rem}[thm]{Remark}

\numberwithin{equation}{section}
\numberwithin{figure}{section}

\def\ge{\geqslant}
\def\le{\leqslant\nobreak}

\def\={\equiv}

\def\R{{\mathbin{\mathbb R}}}

\def\al{\alpha}

\def\de{\delta}

\def\ep{\epsilon}

\def\Si{\Sigma}

\def\d{\partial}

\def\-{\setminus}

\def\ov{\overline}

\begin{document}
\title{A Uniqueness Theorem for Gluing Calibrated Submanifolds}
\author{Yohsuke Imagi}
\date{Kavli IPMU (WPI)}
\maketitle
\section*{Abstract}
`Gluing' is a technique of constructing solutions to non-linear (elliptic) partial differential equations such as Yang--Mills equations, minimal surface equations and Einstein equations. Calibrated submanifolds are a certain class of minimal surfaces, and there are various examples of them constructed by the gluing technique. We have existence theorems in that sense, but there seems to have been no uniqueness theory for higher-dimensional ones such as special Lagrangian submanifolds, which we discuss in the present paper.

\section{Introduction}\l{introduction}
As we have mentioned above `gluing' is a technique of constructing solutions to non-linear (elliptic) partial differential equations such as Yang--Mills equations, minimal surface equations and Einstein equations.
Solutions constructed by the gluing technique are usually parametrized by small $s>0$ and tending to something {\it singular} as $s\to+0$; for example Taubes \cite{Taubes} constructed a one-parameter family of Yang--Mills ASD (anti-self-dual) instantons $A_s$ with curvature tending to a $\de$-function as $s\to+0$. There are many other examples of Yang--Mills instantons, minimal surfaces and Einstein metrics constructed by the gluing technique, including calibrated submanifolds (which are a certain class of minimal surfaces);
for instance various authors \cite{Butscher,Joyce,D. Lee, Y. Lee} constucted  various kinds of special Lagrangian submanifolds (which are a higher-dimensional example of calibrated submanifolds).

What we shall study in the present paper is a {\it uniqueness} problem: given a singular solution and a family of (non-singular) solutions parametrized by $s>0$ and tending to the singular one as $s\to+0$ then need they be re-constructed by the gluing technique?

The answer is `yes' in the situation of Taubes: all ASD instantons with curvature close to a $\de$-function may be re-constructed by the method of Taubes, which was proved by Donaldson \cite{Donaldson}.
Something similar holds for pseudo-holomorphic curves in symplectic manifolds, and they give a key step to the definition of Donaldson invariants (which `count' ASD instantons), Gromov--Witten invariants (which `count' pseudo-holomorphic curves) and Floer homologies in Yang--Mills gauge theory or in symplectic geometry.

There seems to have been no such kind of uniqueness results proved for calibrated submanifolds of higher dimension; pseudo-holomorphic curves are calibrated submanifolds of dimension $2$ and by `higher' we mean $\ge3$.

We shall now recall an outline of the proof of Donaldson.
Let $A_s$ be an instanton whose curvature is close to a $\de$-function supported at a point $x$ in a manifold $X$ (of dimension $4$ and supposed to be compact). There are mainly three things to do:
\i
\item[(i)]
One first proves that $A_s$ tends to the trivial instanton over each compact subset of $X\-\{x\}$ and that there exists $\ep_s>0$ such that if $A_s$ is re-scaled by $\ep_s^{-1}$ about $x$ then the re-scaled instanton $\ep_s^{-1}A_s$ will tend to an instanton $B$ over $T_xX\cong\R^4$ decaying to the trivial instanton at infinity.
\item[(ii)]
There is a well-known classification result for such instantons, which implies that $B$ is a basic instanton used by Taubes (which is unique up to re-scaling). Thus $A_s$ will be close to the trivial instanton on $M\-U$ for some neighbourhood $U$ of $x$ in $X$ and to the re-scaled instanton $\ep_sB$ on a smaller neighbourhood $U_s$ of $x$ in $U$, but we have not seen yet the behaviour of $A_s$ in $U\-U_s$.
\item[(iii)]
We may suppose that $U$ and $U_s$ are open balls about $x$ in $X$.
Since $B$ decays at infinity in $\R^4$ it follows that $\ep_sB$ is close to the trivial instanton near the boundary of $U_s$ and so $A_s$ is close to the trivial instanton near the two boundaries of the annulus $U\-U_s$.
The final step is to prove that $A_s$ is close to the trivial instanton over the whole annulus $U\-U_s$, which will readily imply that $A_s$ is gauge-equivalent to one of the instantons constructed by the gluing technique.
\i

We wish to recover the steps (i)--(iii) for calibrated submanifolds in place of instantons, but it seems too difficult to do in general. We shall therefore focus upon a situation of Joyce \cite{Joyce}. Let $L$ be a compact special Lagrangian submanifold with isolated conical singularities in his sense.
Joyce proved that if there are local models of desingularizing the tangent cones to $L$ then one can glue them to $L$ to get compact special Lagrangian submanifolds (without singularities) parametrized by $s>0$ and tending to $L$ as $s\to+0$ (in the sense of geometric measure theory).

We wish to prove that all compact special Lagrangian submanifolds tending to $L$ may be re-constructed by the method of Joyce. For that purpose it suffices to prove analogues to (i)--(iii) above, and in the present paper we shall prove the analogue to (iii). It seems difficult to prove analogues to (i) and (ii) in full generality, but is doable in some interesting situations, which we do in a sequel to the present paper \cite{I}
(the analogue to (i) holds for special Lagrangian Jacobi-integrable smooth cones and the analogue to (ii) hold for stable $T^2$-cones, which are automatically special Lagrangian Jacobi-integrable).

The analogue to (iii) may be stated as follows:
\begin{thm}
\l{main result}
Let $B(\rho)$ be the open ball of radius $\rho>0$ about $\bm0$ in $\R^n$.
Let $s\in(0,1)$, let $\phi$ be a calibration of degree $m$ on $\R^n$,
and let $CX$ be a $\phi$-calibrated smooth cone in $\R^n$ with $X\=CX\cap S^{n-1}$ being a compact submanifold of $S^{n-1}$. Let $M$ be a properly-embedded $\phi$-calibrated submanifold of $B(1)\-\,\ov{\!B(s)\!}\,$ with $\d M$ being a smooth hypersurface of $\,\ov{\!M}$ contained in $\d B(1)\cup\d B(s)$. Let $\d M\cap\d B(1)$ and $\d M\cap\d B(s)$ be $C^1$-close to $CX\cap\d B(1)$ and $CX\cap\d B(s)$ respectively. Then $M$ is $C^1$-close to $CX\cap B(1)\-\,\ov{\!B(s)\!}\,$.
\end{thm}
\begin{rem}
This holds for general calibrated submanifolds which need not be special Lagrangian. The cone $CX$ plays the r\^ole of the trivial instanton in the step (iii) above, and $M$ plays the r\^ole of $A_s$.
\end{rem}
\begin{rem}
We do not suppose anything particular about the behaviour of $M$ away from $\d M$, but the conclusion of Theorem \ref{main result} implies that $M$ is diffeomorphic to $CX\cap B(1)\-\,\ov{\!B(s)\!}\,$.
\end{rem}
\begin{rem}
The statement above will be refined in Theorem \ref{key step} below.
We have supposed so far that the metric is flat and the calibration is constant, but shall deal with more general metrics and calibrations in Theorem \ref{key step}; it will be necessary for the situation of gluing special Lagrangian submanifolds in Calabi--Yau manifolds where the metric need not be flat and the calibration need not be constant.
In Theorem \ref{key step} we shall also say how close $M$ is to $CX$.
\end{rem}
The proof of Theorem \ref{main result} may be sketched as follows.
Donaldson used a method of Uhlenbeck \cite{Uhlenbeck} in the step (iii) above,
and we shall use a method of Simon \cite{Simon,Simon2} for the proof of Theorem \ref{main result};
Uhlenbeck proved a removable singularity theorem for Yang--Mills instantons in dimension $4$ and Simon proved the uniqueness of multiplicity $1$ smooth tangent cones to minimal surfaces.
Let $M$ be a minimal surface of dimension $m$ in $\R^n$ with an isolated singularity at $\bm0$. It is well-known that $\area\bigl(M\cap B(\rho)\bigr)/\rho^m$ is a monotone non-decreasing function in $\rho$
which plays a central r\^ole in the proof of Simon.
In the situation of Theorem \ref{main result} however we have to work in {\it annuli} instead of balls. We shall therefore make the following version of monotonicity formula.

Let $r$ be the radius function on $\R^n$. For each compact $(m-1)$-dimensional submanifold $\Si$ of $S^{n-1}$ let
\e
F(\Si)\=\int_\Si r^{1-m}\frac{\d}{\d r}\lrcorner\phi.\e
Let $M$ be as in Theorem \ref{main result}.
Then $F\bigl(M\cap\d B(\rho)\bigr)$ will be a monotone non-decreasing in $\rho$ (we note that for $\rho$ generic $M\cap\d B(\rho)$ is a submanifold of $\d B(\rho)\cong S^{n-1}$ which makes $F\bigl(M\cap\d B(\rho)\bigr)$ well-defined almost everywhere in $\rho$);
the functional $F$ is a higher-dimensional analogue of Hofer's functional for pseudo-holomorphic curves in symplectizations of contact manifolds~\cite[pp534--539]{Hofer}. We also note that $F$ is similar to the Chern--Simons functional in the step (iii) above.

Morally speaking $M\cap\d B(\rho)$ behaves like a gradient flow of $F$ where $\rho$ may be regarded as `time'. Theorem \ref{main result} assumes that the flow starts and ends near $X$, and concludes that the flow stays near $X$ for all time. It will follow from Simon's estimates including a version of \L{}ojasiewicz inequality \cite[Lemma~1, p542]{Simon}.

The remainder of the paper will be organized as follows:
\i
\item
In \S\ref{statement section} we state the refined version of Theorem \ref{main result}.
\item
In \S\ref{lemma section 1} we prove the monotonicity formula for $F_\rho$. We shall have error terms in general if the metric is non-flat or the calibration is non-constant.
\item
In \S\ref{lemma section 2} we show how to use Simon's estimates \cite{Simon,Simon2} including a version of \L{}ojasiewicz inequality \cite[Lemma~1, p542]{Simon}.
\item
In \S\ref{proof section}
we complete the proof of Theorem \ref{main result}.
\i
\section*{Acknowledgements}
I was a PhD student at Kyoto University and supported by Grant-in-Aid for JSPS fellows~(22-699) whilst writing the original version of the present paper.
I would like to thank my supervisor Kenji Fukaya for useful conversations.
\section{Statement of Main Result}
\label{statement section}

We begin with a review of calibrated geometry~\cite{Harvey and Lawson}.
Let $W$ be a Riemannian manifold.
An $m$-form $\phi$ on $W$ is said to be of comass~$\le1$ if
$\phi (v_1,\dots,v_m)\le 1$
for every orthonormal vector fields $v_1,\dots,v_m$ on $W$.
A closed $m$-form of comass~$\le1$ on $W$ is called a calibration of degree $m$ on $W$.
Let $\phi$ be a calibration of degree $m$ on $W$.
Let $M$ be an oriented submanifold of $W$.
We call $M$ a $\phi$-submanifold of $W$ if $\phi|_M$ is the volume form of $M$.
By a theorem of Harvey and Lawson~\cite{Harvey and Lawson},
$\phi$-submanifolds of $W$ are minimal submanifolds of $W$.

We shall set up the notation which we use in the statement of Theorem~\ref{key step} below.
Let $g'$ be the Euclidean metric on $\mathbb{R}^n$, i.e.,
\[g'=dy^1\otimes dy^1+\dots+dy^n\otimes dy^n\]
in the coordinates $(y^1,\dots,y^n)$ on $\mathbb{R}^n$.
Let $\phi'$ be a calibration of degree $m$ on $(\mathbb{R}^n,g')$.
Suppose $\phi'$ is parallel, i.e.,
\[\phi'=\phi'_{i_1\dots i_m}dy^{i_1}\wedge\dots\wedge dy^{i_m}\]
for some $\phi'_{i_1\dots i_m}\in\mathbb{R}$.
Let $r$ be the radial coordinate $|\bullet|$ on $(\mathbb{R}^n\setminus\{0\},g')$. 
Set
\begin{equation}
\label{psi'}
\psi'=\left( \partial_r \lrcorner \phi' \right)|_{S^{n-1}},
\end{equation}
where $\partial_r$ is the vector field $\partial/\partial r$,
$\lrcorner$ is the interior product of vector fields with differential forms,
and $S^{n-1}$ is the unit sphere of $(\mathbb{R}^n,g')$.
For every orthonormal vector fields $v_1,\dots,v_{m-1}$ on $S^{n-1}$, we have
\begin{equation}
\label{psi and phi}
\psi'(v_1,\dots,v_{m-1})=\phi'(\partial_r,v_1,\dots,v_{m-1}) \le 1
\end{equation}
since $\partial_r,v_1,\dots,v_{m-1}$ are orthonormal. Therefore,
$\psi'$ is an $(m-1)$-form of comass~$\le 1$ on $S^{n-1}$.
Let $X$ be a oriented submanifold of $S^{n-1}$.
We call $X$ a $\psi'$-submanifold if $\psi'|_X$ is the volume form of $X$.
%
%
\begin{prop}
\label{psi minimal}
$\psi'$-submanifolds of $S^{n-1}$ are minimal submanifolds of $S^{n-1}$.
\end{prop}
\begin{proof}
Let $X$ be a $\psi'$-submanifold of $S^{n-1}$.
Set \[CX=\{rx\in\mathbb{R}^n|r\in (0,\infty),x\in X\}.\]
Then, by \eqref{psi and phi}, $CX$ is a $\phi'$-submanifold of $(\mathbb{R}^n,g')$.
Therefore, $CX$ is a minimal submanifold of $(\mathbb{R}^n,g')$.
Therefore, $X$ is a minimal submanifold of $S^{n-1}$.
\end{proof}
Let $I$ be an open interval of $(0,\infty)$, and $X$ a submanifold of $S^{n-1}$.
We embed $I\times S^{n-1}$ into $\mathbb{R}^n$ by $(r,y)\mapsto ry$.
Let $\nu$ be a normal vector field on $I\times X$ in $(I\times S^{n-1},g')$.
Set
\begin{align*}
\|\nu\|_{C^0_\mathrm{cyl}}=\sup_{I\times X} |\nu|/r ,\;
\|\nu\|_{C^1_\mathrm{cyl}}=\sup_{I\times X} \bigl( |\nu|/r+|D\nu| \bigr),
\end{align*}
where
$D\nu$ is the covariant derivative of $\nu$.
These are induced by the cylindrical metric $g'/r^2$ on $(0,\infty)\times S^{n-1}$.
Set
\begin{equation*}
G_{\mathrm{cyl}}(\nu)
=\Bigl\{ \frac{r}{\sqrt{r^2+|\nu(rx)|^2}}\bigl(rx+\nu(rx)\bigr) \Bigm| r\in I, x\in X \Bigr\}.
\end{equation*}

We are ready now to refine the statement of Theorem \ref{main result}:
\begin{thm}
\label{key step}
Let $B(\rho)$ be the ball of radius $\rho$ about $\bm0$ in $(\R^n,g')$.
Let $\phi'$ be a parallel calibration of degree $m$ on the Euclidean space $(\mathbb{R}^n,g')$,
and $\psi'$ the $(m-1)$-form \eqref{psi'} on the unit sphere $S^{n-1}$ of $(\mathbb{R}^n,g')$.
Let $X$ be a compact $\psi'$-submanifold of $S^{n-1}$.
Let $0<l<1$.
Then, there exist $\epsilon_0, C_0,c_0>0$ depending only on $l,m,n,X,\phi'$ such that if:
\begin{itemize}
\item[$\mathrm{(A0)}$]
$0<\epsilon<\epsilon_0;$
\item[$\mathrm{(A1)}$]
$0<a_0<b_0<a_1<b_1$, $a_0/b_0=a_1/b_1=l;$
\item[$\mathrm{(A2)}$]
$g$ is a Riemannian metric on $B(b_1)$ with
\[\|g-g'\|_{C^1(B(b_1))}\le \epsilon,\;\|g-g'\|_{C^2(B(b_1))}\le1\]
with respect to $g';$
\item[$\mathrm{(A3)}$]
$\phi$ is a calibration on $(B(b_1),g)$ with
\[(1+\log{\frac{b_1}{a_0}})\sup_{B(b_1)}|\phi -\phi'|\le \epsilon,\]
where $|\bullet|$ is with respect to $g';$
\item[$\mathrm{(A4)}$]
$M$ is a closed subset of $(a_0, b_1)\times S^{n-1}$, and
$M$ is a $\phi$-submanifold with respect to $g;$
\item[$\mathrm{(A5)}$]
there exists a normal vector field $\nu_i$
on $(a_i,b_i)\times X$ in $((a_i,b_i)\times S^{n-1},g'/r^2)$, where $i=0,1,$ such that
\begin{align*}
& M\cap ((a_i, b_i)\times S^{n-1})=G_{\mathrm{cyl}}(\nu_i)
\text{ with }\| \nu _i \|_{C^1_\mathrm{cyl}}\le \epsilon,
\end{align*}
\end{itemize}
then there exists a normal vector field $\nu$ on $(a_0,b_1)\times X$ in $((a_0,b_1)\times
S^{n-1},g'/r^2)$ such that
\begin{equation}
\label{conclusion of key step}
M= G_{\mathrm{cyl}}(\nu)
\text{ with }\| \nu \|_{C^1_\mathrm{cyl}}\le
C_0\epsilon^{c_0}.
\end{equation}
\end{thm}
\begin{rem}
One sufficient condition for (A3) to hold is that we have $\phi|_{\bm0}=\phi'$ and $a_0=s^\al,\,b_1=s$ for some $s>0$ small enough and $\al\in(0,1)$ independent of $s$; if so we have
\[(1+\log{\frac{b_1}{a_0}})\sup_{B(b_1)}|\phi -\phi'|=(1+(1-\al)\log s)O(s)\]
which tends to $0$ as $s\to+0$.
\end{rem}
\section{A Monotonicity Formula}
\label{lemma section 1}
In this section we prove a monotonicity formula for calibrated submanifolds of annuli; see Proposition~\ref{dpsi-energy proposition}.
This is a higher-dimensional analogue of an energy estimate of Hofer~\cite[pp534--539]{Hofer} for pseudo-holomorphic curves in symplectizations of contact manifolds.

Let $g$ be a Riemannian metric on $\mathbb{R}^n$,
and $\phi$ a calibration of degree $m$ on $(\mathbb{R}^n,g)$.
%
\begin{prop}\label{normal to calibrated submanifold}
Let $M$ be a $\phi$-submanifold of $(\mathbb{R}^n,g)$.
If $\nu$ is a normal vector field on $M$ in $(\mathbb{R}^n,g)$, then we have
\[(\nu\lrcorner\phi)|_M=0.\]
\end{prop}
\begin{proof}
It suffices to prove that
for every point $p\in M$ and orthonormal vectors $v_1, \dots, v_{m-1}\in T_pM$, we have
\begin{equation}
\label{normal zero}
\phi_p(\nu_p,v_1,\dots,v_{m-1})=0.
\end{equation}
Choose $v\in T_pM$ so that $\phi_p(v, v_1, \dots, v_{m-1})=1$.
Consider
\begin{equation*}
\label{normal phi}
t\mapsto \phi_p((\sin t)\nu_p+(\cos t)v,v_1,\dots,v_{m-1}).
\end{equation*}
By the definition of calibration, this attains maximum $1$ at $t=0$.
Differentiating it at $t=0$, we have \eqref{normal zero}.
\end{proof}
Let $g'$ be the Euclidean metric on $\mathbb{R}^n$.
Let $r$ be the radial coordinate on the Euclidean space $(\mathbb{R}^n,g')$, and $\partial_r$ the vector field $\partial/\partial_r$.
In the same way as Harvey and Lawson~\cite[Lemma~5.11, II.5]{Harvey and Lawson}, we shall prove the following
%
\begin{prop}
\label{energy}
Let $M$ be a $\phi$-submanifold of $(\mathbb{R}^n,g)$.
Then, we have
\begin{equation}
\label{energy formula}
\langle\overrightarrow{TM}, \partial_r\lrcorner dr \wedge\phi\rangle=|\mathrm{pr}_{TM^{\perp}}\partial_r|^2,
\end{equation}
where $\langle\bullet,\bullet\rangle$ is the canonical pairing of poly-vector fields and differential forms,
$\overrightarrow{TM}$ is the $m$-vector field on $M$ dual to $\phi|_M$,
$r=|\bullet|$ is with respect to the Euclidean metric $g'$,
and $\mathrm{pr}_{TM^{\perp}}$ is the projection of $\mathbb{R}^n$ onto the normal bundle
of $M$ in $(\mathbb{R}^n,g)$.
\end{prop}
\begin{proof}
By Proposition~\ref{normal to calibrated submanifold}, we have
\begin{align*}
\langle\nu\wedge\overrightarrow{TM},dr\wedge\phi\rangle
=\langle\nu, dr\rangle\langle\overrightarrow{TM},\phi\rangle,
\text{ where }\nu=\mathrm{pr}_{TM^{\perp}}\partial_r.
\end{align*}
This proves \eqref{energy formula}.
\end{proof}
Set
\begin{equation}\label{psi}
\psi=\frac{m}{r^m}\int_0^r(\partial_r\lrcorner\phi)dr.
\end{equation}
\begin{prop}\label{phi=d}
$\psi$ is an $(m-1)$-form on $\mathbb{R}^n\setminus\{0\}$ such that
\begin{equation}
\label{phi=d formula}
\phi=d\left( \frac{r^m}{m}\psi \right).
\end{equation}
\end{prop}
\begin{proof}
Set $\chi=\partial_r\lrcorner\phi$, and $\omega=\partial_r \lrcorner dr \wedge\phi$.
Then, we have
\begin{equation}\label{decomposition of phi}
\phi=dr\wedge\chi+\omega.
\end{equation}
Since $\partial_r\lrcorner\chi=\partial_r\lrcorner\omega=0$, we may regard $\chi$ and $\omega$ as smooth families of
differential forms on $S^{n-1}$.
By the definition of calibration, $d\phi=0$.
Therefore, we have
\begin{equation}
\label{dchi}
d_{S^{n-1}}\chi=\partial_r\omega,
\end{equation}
where $d_{S^{n-1}}$ is the exterior differentiation on $S^{n-1}$.
By \eqref{decomposition of phi} and \eqref{dchi}, we have
\[ \phi=d\left( \int_0^r \chi dr \right)=d\left( \int_0^r (\partial_r\lrcorner\phi) dr \right).\]
By \eqref{psi}, this proves \eqref{phi=d formula}.
\end{proof}
%
%
%
Let $\phi'$ a parallel calibration of degree $m$ on the Euclidean space $(\mathbb{R}^n,g')$,
Set
\begin{equation}\label{psi'2}
\psi'=r^{1-m}\partial_r\lrcorner\phi'.
\end{equation}
Then, \eqref{psi} holds with $\phi',\psi'$ in place of $\phi,\psi$ respectively.
%
%

We shall prove a monotonicity formula with an error term.
When $\phi=\phi'$,
it has no error term.
\begin{prop}\label{dpsi-energy proposition}
There exists $C_{m,n}>0$ depending only on $m,n$ such that 
\begin{equation}
\label{dpsi-energy}
 \left| m^{-1}d\psi  - r^{-m}\partial_r\lrcorner dr\wedge\phi \right|_{\mathrm{cyl}}
\le C_{m,n}\sup|\phi-\phi'|,
\end{equation}
where $|\bullet|_{\mathrm{cyl}}$ is with respect to the metric $g'/r^2$.
\end{prop}
\begin{proof}
By \eqref{phi=d formula} and \eqref{psi'2}, we have
\begin{equation}
\label{dpsi-energy'}
m^{-1}d\psi  - r^{-m}\partial_r\lrcorner dr\wedge\phi
=dr/r\wedge (r^{1-m}\partial_r\lrcorner\phi-r^{1-m}\partial_r\lrcorner\phi'+\psi'-\psi).
\end{equation}
By \eqref{psi} and \eqref{psi'2}, we have
\begin{align*}
&|r^{1-m}\partial_r\lrcorner\phi-r^{1-m}\partial_r\lrcorner\phi'|_{\mathrm{cyl}}\le c\sup|\phi-\phi'|,\\
&|\psi-\psi'|_{\mathrm{cyl}}\le c\sup|\phi-\phi'|
\end{align*}
for some $c>0$ depending only on $m,n.$
Therefore, by \eqref{dpsi-energy'}, we have \eqref{dpsi-energy}.
\end{proof}
%
%
%
%
We shall prove a proposition which we use in the proof of Lemma~\ref{first lemma} below.
We also use it in the key step to proof of the main result of this paper.
\begin{prop}
\label{vol estimate}
Let $M$ be a $\phi$-submanifold of $(\mathbb{R}^n,g)$,
and suppose $M$ is a closed subset of $(a,b)\times S^{n-1}$, where $(a,b)\times S^{n-1}$ is embedded into $\mathbb{R}^n$ by $(r,y)\mapsto ry$.
There exist $\epsilon_{m,n},C'_{m,n}>0$ depending only on $m,n$ such that
if
\begin{align}
\label{epsilonlog}
(1+m\log{\frac{b}{a}})\sup_{(a,b)\times S^{n-1}}|\phi-\phi'|\le \epsilon_{m,n},\;
\sup_{(a,b)\times S^{n-1}}|g-g'|\le 1,
\end{align}
then we have
\begin{equation}
\label{vol estimate 2}
\begin{split}
\Vol(M,g/{r^2})
&\le
C'_{m,n}\log{\frac{b}{a}}\limsup_{r\to b}\left|\int_{M\cap \{r\}\times S^{n-1}}\psi\right|\\
&+C'_{m,n}(1+m\log{\frac{b}{a}})\int_M |\mathrm{pr}_{TM^{\perp}}\partial_r |^2\dVol(M,g/r^2).
\end{split}
\end{equation}
\end{prop}
\begin{proof}
By \eqref{phi=d formula}, we have
\begin{equation*}
\label{vol=int}
\Vol(M,g/r^2)
=\int_M\phi/r^m
=\int_M (dr/r)\wedge \psi+{m}^{-1}\int_M d\psi.
\end{equation*}
By \eqref{dpsi-energy}, we have
\begin{equation*}
\label{int dpsi}
m^{-1}\int_M d\psi\le \int_M |\mathrm{pr}_{TM^{\perp}}\partial_r|^2\dVol(M,g/r^2)
+C_{m,n}\sup|\phi-\phi'|\Vol(M,g'/r^2).
\end{equation*}
By \eqref{dpsi-energy} and \eqref{energy formula}, we have
\begin{equation*}
\label{Psi proposition}
\begin{split}
&\int_M (dr/r)\wedge \psi
\le \log{\frac{b}{a}} \limsup_{r\to b}\left|\int_{M\cap \{r\}\times S^{n-1}}\psi
+\int_{M\cap ([a,r]\times S^{n-1})}d\psi\right|\\
&\le m\log{\frac{b}{a}}  \limsup_{r\to b} \left| \int_{M\cap \{r\}\times S^{n-1}}m^{-1}\psi\right|
+\int_{M}|\mathrm{pr}_{TM^{\perp}}\partial_r|^2\dVol(M,g'/r^2)\\
&+mC_{m,n}\log{\frac{b}{a}}\sup|\phi-\phi'|\Vol(M,g'/r^2).
\end{split}\end{equation*}
Thus, we have
\begin{equation*}
\label{vol estimate 3}
\begin{split}
&\Vol(M,g/{r^2})\\
&\le m\log{\frac{b}{a}}\limsup_{r\to b}\left|\int_{M\cap\{r\}\times S^{n-1}}\psi\right|
+(1+m\log{\frac{b}{a}})\int_M |\mathrm{pr}_{TM^{\perp}}\partial_r|^2\dVol(M,g/r^2)\\
&+C_{m,n}(1+m\log{\frac{b}{a}}) \sup |\phi-\phi'|\Vol(M,g'/r^2).
\end{split}
\end{equation*}
By \eqref{epsilonlog}, we have
\begin{align*}
C_{m,n}(1+m\log{\frac{b}{a}}) \sup |\phi-\phi'|\Vol(M,g'/r^2)\le (1/2)\Vol(M,g/r^2).
\end{align*}
Thus, we have \eqref{vol estimate 2}.
\end{proof}
We shall prove a lemma which we use in the key step to the proof of the main result of this paper.
It is similar to a lemma of Simon~\cite[Lemma~3, p561]{Simon}.
We however use the monotonicity formula for $\phi$-submanifolds of annuli.
\begin{lem}
\label{first lemma}
Let $\phi'$ be a parallel calibration of degree $m$ on the Euclidean space $(\mathbb{R}^n,g')$,
and let $\psi'$ be as in \eqref{psi'2}.
Let $X$ be a compact $\psi'$-submanifold of ${S}^{n-1}$.
Let $\epsilon>0$, and $0<\lambda<\lambda''<\lambda'<1$.
Then, there exists $\delta>0$ such that
if:
\begin{itemize}
\item[\rm(P1)]
$g$ is a Riemannian metric on $B(1)$ with $\|g-g'\|_{C^1(B(1))}\le \delta$,
where $\|\bullet\|_{C^1}$ is with respect to $g'$, and $B(1)$ is the unit ball of $(\mathbb{R}^n,g');$
\item[\rm(P2)]
$\phi$ is a calibration on $(B(1),g)$ with $\sup_{B(1)}|\phi -\phi'|\le \delta$,
where $|\bullet|$ is with respect to $g';$
\item[\rm(P3)]
$M$ is a $\phi$-submanifold of $(\mathbb{R}^n,g)$, and
$M$ is a closed subset of $(\lambda,1)\times S^{n-1}$, where $(\lambda,1)\times S^{n-1}$ is embedded into $\mathbb{R}^n$ by $(r,y)\mapsto ry;$
\item[\rm(P4)]
there exists a normal vector field $\nu$
on $(\lambda',1)\times X$ in $((\lambda',1)\times S^{n-1},g'/r^2)$ such that
\begin{equation*}
M\cap ((\lambda', 1)\times S^{n-1})=G_{\mathrm{cyl}}(\nu)\text{ with }\| \nu \|_{C^1_\mathrm{cyl}}\le \delta
\end{equation*}
in the notation of Section~$\ref{statement section};$
\item[\rm(P5)]
$\int_M |\mathrm{pr}_{TM^{\perp}}\partial_r|^2\dVol(M,g/r^2) \le \delta,$
\end{itemize}
then there exists a normal vector field $\nu'$ on $(\lambda,1)\times S^{n-1}$
in $((\lambda,1)\times S^{n-1},g'/r^2)$ such that
\begin{equation*}
M= G_{\mathrm{cyl}}(\nu')\text{ with }
\| \nu '|_{(\lambda'',\lambda')\times S^{n-1}} \|_{C^{1,1/2}_{\mathrm{cyl}}} \le \epsilon,
\end{equation*}
where $C^{1,1/2}_{\mathrm{cyl}}$ is the H\"older space with respect to the
metric $g'/r^2$ on $(\lambda,1)\times S^{n-1}$.
\end{lem}
\begin{proof}
Suppose there does not exist such $\delta$.
Then, for every $j=2,3,4,\dots,$ there exist $g_j,\phi_j,M_j$ such that
(P1), (P2), (P3), (P4) and (P5) hold with $\delta=1/j$, and the following holds: 
\begin{itemize}
\item[(P6)]
there does not exist any normal vector field $\nu_j'$ on
$(\lambda'',\lambda')\times X$ in $((\lambda'',\lambda')\times S^{n-1},g'/r^2)$ such that
\begin{align*}
& M_j= G_{\mathrm{cyl}}(\nu_j') \text{ with }\| \nu _j' \|_{C^{1,1/2}_\mathrm{cyl}}\le \epsilon.
\end{align*}
\end{itemize}
By (P1), (P2) and (P3),
we may apply Proposition~\ref{vol estimate}.
Therefore, by \eqref{vol estimate 2}, (P4) and (P5), we have
\begin{equation}
\label{boundedness}
 \sup_{j=2,3,4,\dots} \Vol(M_j,g_j/r^2)<\infty.
\end{equation}
%
%
%
Therefore, by (P1), we have
\begin{equation*}
\label{boundedness'}
 \sup_{j=2,3,4,\dots} \Vol(M_j,g')<\infty.
\end{equation*}
By (P1) and (P3), we have
\begin{equation}
\label{mean curvature}
\lim_{j\to\infty}\left(\text{the mean curvature of }M_j\text{ in }((\lambda,1)\times {S}^{n-1},g')\right)=0 
\end{equation}
in the $C^0$-topology.
Thus, by
Allard's compactness theorem~\cite[Theorem 5.6]{Allard},
there exists a subsequence
$M_{j_k}$ converging as varifolds to some rectifiable varifold $M_{\infty}$ in
$((\lambda,1)\times{S}^{n-1},g')$.

Let $\|M_{\infty}\|$ be the Radon measure on $((\lambda,1)\times{S}^{n-1},g')$ induced by $M_{\infty}$.
%
%
%
We shall prove
\begin{equation}
\label{scaling invariance}
a^m\|M_{\infty}\|(a^{-1}E)=\|M_{\infty}\|(E)
\end{equation}
for every $a>0,E\subset(\lambda,1)\times{S}^{n-1}$ with $aE\subset (\lambda,1)\times{S}^{n-1}$.
It suffices to prove
\begin{equation}
\label{approximation by smooth function}
\frac{d}{da}a^{m}\int_{(\lambda,1)\times{S}^{n-1}}f(ar)h d\|M_{\infty}\| =0
\end{equation}
for every smooth functions $h:{S}^{n-1}\rightarrow[0,\infty)$ and $f:(\lambda,1)\rightarrow[0,\infty)$ with $a(\mathrm{supp}{f})\subset (\lambda,1)$.
By \eqref{dpsi-energy}, (P2), (P5) and \eqref{boundedness}, we have
\[\lim_{j\to\infty}\int_{M_j}d\psi_j\to0,\]
where
$\psi_j$ is as in \eqref{psi} with $\phi_j$ in place of $\phi$.
Therefore, by (P3) and \eqref{phi=d formula}, we have
\begin{equation*}
\label{11111111}
\begin{split}
\text{the left-hand side of \eqref{approximation by smooth function}}
&=\frac{d}{da}\lim_{k\to \infty}
a^{m}\int_{M_{j_k}}f(ar)h d(\frac{r^m}{m}\psi_{j_k})\\
&=\lim_{k\to \infty}
\int_{M_{j_k}}\frac{d}{da}((ar)^m f(ar))\frac{dr}{r}\wedge h\psi_{j_k}\\
&=\lim_{k\to \infty}
\int_{M_{j_k}}a^{-1}\frac{d}{dr}((ar)^{m}f(ar)) dr\wedge h\psi_{j_k} \\
&=\lim_{k\to \infty}
-\int_{M_{j_k}}a^{-1}(ar)^{m}f(ar) dh\wedge \psi_{j_k}.
\end{split}
\end{equation*}
Therefore, by \eqref{psi'2}, (P2) and \eqref{boundedness}, we have
\[\text{the left-hand side of \eqref{approximation by smooth function}}
=\lim_{k\to \infty}
-\int_{M_{j_k}}a^{-1}{(ar)^{m}f(ar)}r^{1-m}\partial_r\lrcorner (dh\wedge \phi_{j_k}).\]
By Proposition~\ref{normal to calibrated submanifold}, we have
\begin{equation*}
\int_{M_{j_k}}r^{-m}\partial_r\lrcorner (dh\wedge \phi_{j_k})
=\int_{M_{j_k}}\langle \mathrm{pr}_{TM_{j_k}^{\perp}}\partial_r, dh \rangle\dVol(M_{j_k},g_{j_k}/r^2).
\end{equation*}
%
This converges to $0$ by (P5) and \eqref{boundedness}.
Thus, we have \eqref{approximation by smooth function}.
This proves \eqref{scaling invariance}.
%

By (P4),
the restriction of $M_{\infty}$ to $(\lambda',1)\times S^{n-1}$
is equal to $(\lambda',1)\times X$ as varifolds in $((\lambda',1)\times S^{n-1},g')$.
Therefore, by \eqref{scaling invariance}, we have
\[M_{\infty}=(\lambda,1)\times X
\text{ as varifolds in }((\lambda,1)\times {S}^{n-1},g').\]
Therefore, $M_{j_k}$ converges to $(\lambda,1)\times X$ as varifolds
in $((\lambda,1)\times {S}^{n-1},g')$.
Therefore, by \eqref{mean curvature} and Allard's regularity theorem~\cite[Theorem 8.19]{Allard},
$M_{j_k}$ converges to $(\lambda,1)\times X$ in the local $C^{1,1/2}$-topology in $(\lambda,1)\times S^{n-1}$.
This contradicts (P6), which completes the proof of Lemma~\ref{first lemma}.
\end{proof}
\section{Simon's Estimates}
\label{lemma section 2}
In this section we show how to use Simon's estimates \cite{Simon,Simon2} including a version of \L{}ojasiewicz inequality \cite[Lemma~1, p542]{Simon}.

Let $X$ be a compact smooth Riemannian manifold, $V$ a smooth real vector bundle on $X$ with a fibre metric and a metric connection. 
Let $C^{\infty}_x$ be the space of smooth sections of $V\rightarrow X$. 
Let $E:C^{\infty}_x\rightarrow \mathbb{R}$ satisfy
\begin{equation}
\label{multiple integral}
 Ev=\int_X F(x,v,D_x v )dx
\end{equation}
for every $v\in C^{\infty}_x$, where $D_xv$ is the covariant derivative of $v$,
and $F=F(x,v,p)$ is a $\mathbb{R}$-valued smooth function of $x\in X$, $v\in V|_x$, $p\in T^*_xX\otimes V|_x$.
Suppose $F$ satisfies the following conditions:
\begin{itemize}
\item[(C1)] $(v,p)\mapsto F(x,v,p)$ is a real-analytic function on the vector space $V|_x \oplus \left( T^*_xX\otimes V|_x \right)$ for every $x\in X$;
\item[(C2)] there exists $c>0$ such that for every $x\in X, \xi\in T_x^*X, v\in V|_x,$ 
\[ \frac{d^2}{dh^2}F(x,0,h^2\xi\otimes v)\Big|_{h=0}>c|\xi|^2|v|^2.\]
\end{itemize}
By (C1), one can use the \L{}ojasiewicz estimate~\cite{Lojasiewicz}. This is important in the proof of a result of Simon; for the statement, see Proposition~\ref{Lojasiewicz-Simon} below.
(C2) is called the Legendre--Hadamard condition.
Let $-\grad E:C^{\infty}_x\rightarrow C^{\infty}_x$ be the Euler--Lagrange operator of $E$, i.e.,
\[ \bigl( \grad E(v),v' \bigr)_{L^2_x}=\frac{d}{dh}E(v+hv')\Big|_{h=0}\]
for every $v,v'\in C^{\infty}_x$, where
\begin{equation}
\label{L^2_x}
(v'',v')_{L^2_x}=\int _X \bigl( v''(x),v'(x)\bigr) dx;
\end{equation}
here $\bigl( v''(x),v'(x)\bigr)$ is the inner product on the fibre $V|_x$ at $x\in X$.
Suppose
\begin{equation}
\label{critical point}
\grad E(0)=0,\text{ where }0\in C^{\infty}_x.
\end{equation}
Let $t_0<t_{\infty}$.
Let $C^{\infty}_{t,x}(t_0,t_{\infty})$ be the space of
all smooth sections $u=u(t,x)$ with $u(t,x)\in V|_x$ for every $(t,x)\in(t_0,t_{\infty})\times X$.
Let $C^{k,\mu}_{t,x}(t_0,t_{\infty})$ be the H\"older spaces with respect to the product metric
on $(t_0,t_{\infty})\times X$.
Set $u(t)=u(t,\bullet) \in C^{\infty}_x$ for every $u=u(t,x)\in C^{\infty}_{t,x}(t_0,t_{\infty})$.

We shall state a result of Simon which we use in the proof of Lemma~\ref{second lemma} below.
\begin{prop}[Simon~{\cite[Lemma~1, p542]{Simon}}]
\label{Lojasiewicz-Simon}
There exist $\delta_0,\theta>0$ depending only on $X$, $V$, $E$
such that if $t_0<t_3<t_4<t_{\infty},u\in C^{\infty}_{t,x}(t_0,t_{\infty}),\delta>0$ and if
\begin{equation}
\label{the assumption of Lojasiewicz-Simon}
\begin{split}
& \|u\|_{C^{2,1/2}_{t,x}(t_3,t_4)} \le \delta_0,\\
& \sup_{t\in[t_3,t_4]}\left(E(0)-E\bigl( u(t)\bigr)\right) \le \delta,\\
& \|\partial_t u(t)+\grad E\bigl(u(t)\bigr) \|_{L^2_x}\le (3/4)\|\partial_t u(t)\|_{L^2_x}\text{ for every }t\in[t_3,t_4],
\end{split}
\end{equation}
then we have
\begin{equation*}
\int_{t_3}^{t_4}\|\partial_t u(t)\|_{L^2_x}dt
\le (4/\theta) \bigl( \bigl|  E\bigl( u(t_3)\bigr)-E(0)  \bigr|^{\theta} + \delta ^{\theta} \bigr).
\end{equation*}
Here, $\|\bullet\|_{L^2_x}$ is with respect to \eqref{L^2_x}.
\end{prop}
%

Consider $u=u(t,x)\in C^{\infty}_{t,x}(t_0,t_{\infty})$ satisfying
\begin{equation}
\label{Simon's equation}
 \partial_t^2u-\partial_tu-\grad E(u)+R(u,\partial_tu,\partial_t^2u) =f
\end{equation}
as in Simon~\cite{Simon}, where
$f\in C^{\infty}_{t,x}(t_0,t_{\infty})$ satisfies
\begin{equation}
\label{solution with 00}
\|\partial_t^k f(t)\|_{C^{2}_{x}}\le C_f e^{-2(t-t_0)}\text{ for every }t\in(t_0,t_{\infty}),k=0,1,2
\end{equation}
for some $C_f>0$, and
$R:C^{\infty}_x\times C^{\infty}_x\times C^{\infty}_x\rightarrow C^{\infty}_x$ satisfies
\begin{equation}
\label{remainder}
\begin{split}
 R( v,v^{(1)},v^{(2)} )=& A(x,v,D_x v,v^{(1)})\cdot D_x^2 v \otimes v^{(1)}\\ & +\sum_{(k,l)=(0,1),(1,1),(0,2)}B_{kl}(x,v,D_x v,v^{(1)})\cdot D_x^{l}v^{(k)}
\end{split}\end{equation}
for every $v,v^{(1)},v^{(2)}\in C^{\infty}_x$, where $A=A(x,v,p,q)$, $B_{kl}=B_{kl}(x,v,p,q)$ are smooth functions of $x\in X$, $v\in V|_x$, $p\in T^*_xX\otimes V|_x$, $q\in V|_x$ with
$A(x,v,p,q)\in \mathrm{Hom}(\bigotimes ^2 T_x^*X\otimes V|_x\otimes V|_x,V|_x)$,
$B_{kl}(x,v,p,q)\in \mathrm{Hom}(\bigotimes ^l T_x^*X\otimes V|_x,V|_x)$ and $B_{kl}(x,0,0,0)=0$ for every $x\in X$, $(k,l)=(1,0),(1,1),(2,0)$.
Then,
for every $C_2'>0$, there exists $\delta_4=\delta_4(X,V,E,R,C_2')>0$ such that
if $\|u\|_{C^{1,1/2}_{t,x}(t_0,t_{\infty})}\le \delta_4$, then we have
\begin{equation}
\label{remainder estimate}
|R(u(t), \partial_t u(t), \partial_t^2 u(t)) |
\le C_2'( |\partial_tu(t)| +|D_x\partial_t u(t)|+|\partial_t^2 u(t)|).
\end{equation}
Let $H:C^{\infty}_x\rightarrow C^{\infty}_x$ be the linearized operator of $\grad E$ at $0\in C^{\infty}_x$.
Then, \eqref{Simon's equation} is of the form
\begin{equation}
\label{the linearized equation}
\partial_t^2u-\partial_tu-Hu=\sum_{0\le k+l\le 2}a_{kl}(x,u,D_xu,\partial_tu)\cdot D_x^l\partial_t^ku +f,
\end{equation}
where $a_{kl}=a_{kl}(x,u,p,q)$ are smooth functions of $x\in X$, $v\in V|_x$, $p\in T^*_xX\otimes V|_x$, $q\in V|_x$ with
$a_{kl}(x,v,p,q)\in \mathrm{Hom}(\bigotimes ^l T_x^*X\otimes V|_x,V|_x)$, $a_{kl}(x,0,0,0)=0$
for every $x\in X,0\le k+l\le2$. Therefore, there exists $\delta_2=\delta_2(X,V,E,R)>0$ such that
%
if $u\in C^{\infty}_{t,x}(t_0,t_{\infty})$
with $\|u\|_{C^{1,1/2}_{t,x}(t_0,t_{\infty})}\le \delta_2$, then we have
\begin{equation}
\label{error1}
\max_{0\le k+l\le2}\|a_{kl}(x,u,D_xu,\partial_tu)\|_{C^{0,1/2}_{t,x}(t_0,t_{\infty})}\le \delta_1,
\end{equation}
where $\delta_1=\delta_1(X,V,E)>0$ is given below.
%
By the Legendre--Hadamard condition~(C2),
$\partial_t^2-\partial_t-H$ is elliptic on $C^{\infty}_{t,x}(t_0,t_{\infty})$.
Therefore, there exists $\delta_1=\delta_1(X,V,E)>0$ such that
if $T>0$, if $w,g\in C^{\infty}_{t,x}(-T/3,T/3)$ and if
\begin{equation}
\label{the linear equation}
\partial_t^2w-\partial_tw-Hw=\sum_{0\le k+l\le2}{b}_{kl}(t,x)\cdot D_x^l\partial_t^k w +g
\end{equation}
with $\max_{0\le k+l\le2}\|b_{kl}\|_{C^{0,1/2}_{t,x}(-T/3,T/3)}\le \delta_1$,
then we have
\begin{equation}
\label{the Schauder estimate}
\|w\|_{C^{2,1/2}_{t,x}(-T/5,T/5)}\le C_1\|w\|_{L^2_{t,x}(-T/4,T/4)}+C_1\|g\|_{C^{0,1/2}_{t,x}(-T/4,T/4)}
\end{equation}
for some $C_1=C_1(X,V,E;T)>0$;
here $L^2_{t,x}(t',t'')$ is with respect to the product metric on $(t',t'')\times X$.
\eqref{the Schauder estimate} is a Schauder estimate for elliptic systems;
see Douglis--Nirenberg~\cite{Douglis and Nirenberg} and Morrey~\cite{Morrey}.
%
%

We shall state a proposition which we use in the proof of Lemma~\ref{second lemma} below.
One can prove it in the same way as a result of Simon; see \cite[Lemma~2, p549]{Simon} or \cite[Lemma~3.3, Part II]{Simon2}.
\begin{prop}
\label{growth theorem}
There exist
$h,T_3,\delta_3 >0$ depending only on $X$, $V$, $E$ such that if $T>T_3$, if
$w,g\in C^{\infty}_{t,x}(0,3T)$ satisfy \eqref{the linear equation} with
$\| {b}_{kl}\|_{C^0_{t,x}(0,3T)} \le \delta_3$,
and if
\[\|g\|_{L^2_{t,x}(0,3T)}\le {\delta_3}^{1/3}\|w\|_{L^2_{t,x}(T,2T)}\text{ with }\|w\|_{L^2_{t,x}(0,3T)}<\infty,\]
then we have
\begin{equation*}
\begin{split}
& \|w\|_{L^2_{t,x} (2T,3T) } \le e^{-h T}\|w\|_{L^2_{t,x} (T,2T)}
\Longrightarrow \|w\|_{L^2_{t,x} (T,2T) } \le e^{-h T}\|w\|_{L^2_{t,x} (0,T)},\\
& \|w\|_{L^2_{t,x}(T,2T)} \ge e^{h T} \|w\|_{L^2_{t,x} (0,T)}
\Longrightarrow \|w\|_{L^2_{t,x}(2T,3T)} \ge e^{h T} \|w\|_{L^2_{t,x} (T,2T)},\\
%
&\|w\|_{L^2_{t,x}(T,2T)}\ge e^{-hT}\|w\|_{L^2_{t,x} (0,T)}
\text{ and }\|w\|_{L^2_{t,x}(2T,3T)}\le e^{hT}\|w\|_{L^2_{t,x}(T,2T)}\\
&\Longrightarrow \| w(t) \|_{L^2_x}\le (3/2) \| w(t') \|_{L^2_x} \text{ for every }t,t'\in (T,2T)\\
&\text{and }\| \partial_tw(t) \|_{L^2_x}\le (1/2) \| w(t) \|_{L^2_x}
\text{ for every }t\in (T,2T).
\end{split}
\end{equation*}
\end{prop}
We shall prove a lemma which we use in the key step to the main result of this paper.
It is similar to a result of Simon~\cite[Theorem~1, p534]{Simon}.
Simon's result is an a-priori estimate on $(0,\infty)\times X$.
We however consider $(t_0,t_{\infty})\times X$ with $(t_0,t_{\infty})$ bounded.
We prove the lemma for completeness.
\begin{lem}
\label{second lemma}
Let $X,V,E,R$ be as above. Let $t_0<t_{\infty}$, and $f\in C^{\infty}_{t,x}(t_0,t_{\infty})$ with \eqref{solution with 00}
for some $C_f>0$.
Then,
there exist $\theta,\delta_*,C_*>0$
depending only on $X,V,E,R,C_f$ such that if
$t_*\in(t_0,t_\infty)$, if $u \in C^{\infty}_{t,x}(t_0,t_*)$ satisfies \eqref{Simon's equation} and if
\begin{align}
\label{solution with 11}
& \|u\|_{C^{1,1/2}_{t,x}(t_0,t_*)}\le \delta_*, \\
\label{solution with 22}
& \limsup_{t \to t_0} \|u(t)\|_{L^2_x}\le \delta, \\
\label{solution with 33}
& \sup_{t\in (t_0,t_*)} \Bigl( E(0)-E\bigl( u(t) \bigr) \Bigr) \le \delta,\\
\label{solution with 44}
& \|\partial_tu\|_{L^2_{t,x}(t_0,t_*)}\le \sqrt{\delta}
\end{align}
for some $0<\delta<\min\{1,\delta_*\}$, then we have
\begin{equation}
\label{conclusion of second lemma}
\sup_{t\in(t_0,t_*)}\| u(t) \|_{L^2_x}<C_*\delta^{\theta}.
\end{equation}
\end{lem}
\begin{proof}
By \eqref{solution with 22}, it suffices to prove
\begin{equation}
\label{a-priori estimate}
\int_{t_0}^{t_*} \| \partial_tu(t) \|_{L^2_x}dt<C_*\delta^{\theta}.
\end{equation}
By the Schwartz inequality and \eqref{solution with 44}, for every $(t',t'')\subset(t_0,t_*)$, we have
\begin{equation}
\label{the first}
\begin{split}
\int_{t'}^{t''} \| \partial_tu(t) \|_{L^2_x}dt
\le \sqrt{t''-t'}\|\partial_t u\|_{L^2_{t,x}(t',t'')}
\le \sqrt{(t''-t')\delta}.
\end{split}
\end{equation}
Let $T>0$ be a sufficiently large constant; in the proof of Lemma~\ref{second lemma} a constant means a
real number depending only on $X,V,E,R,C_f$.
If $t_*-t_0<8T$, then by \eqref{the first}, we have \eqref{a-priori estimate};
we may therefore assume $t_*-t_0\ge 8T$.
Choose $t_1, t_6\in(t_0,t_*)$ so that $T\le t_1-t_0\le 2T$, $T\le t_*-t_6\le 2T$ and 
$t_6-t_1=jT$
for some integer $j\ge 4$.
Then, by \eqref{the first}, we have
\begin{align}
\label{t_0 to t_1} &\int_{t_0}^{t_1}\|\partial_t u(t)\|_{L^2_{t,x}}\le \sqrt{T\delta} ,
\; \int_{t_6}^{t_*}\|\partial_t u(t)\|_{L^2_{t,x}}\le \sqrt{T\delta}.
\end{align}
By \eqref{solution with 11}, $u$ satisfies \eqref{the linearized equation} with \eqref{error1}.
Therefore, $u$ satisfies the Schauder estimate~\eqref{the Schauder estimate}. 
Therefore, by \eqref{solution with 11} and \eqref{solution with 00}, we have
\[\|u\|_{C^{2,1/2}_{t,x}(t_1,t_6)}\le C_1'\delta_*+ C_1''e^{-2T}\]
for some constants $C_1',C_1''>0$.
We may therefore assume that
\begin{equation}
\label{solution with 22'}
 \|u\|_{C^{2,1/2}_{t,x}(t_1,t_6)}
 \text{ is sufficiently small}.
 \end{equation}

Differentiating \eqref{Simon's equation} with respect to $t$ and using \eqref{solution with 22'}, we have:
\begin{equation}
\label{error2}
w=\partial_tu,g=\partial_tf\text{ satisfy \eqref{the linear equation} with }
\|{b}_{kl}\|_{C^{0,1/2}_{t,x}(t_1,t_6)}
\text{ sufficiently small}.
\end{equation}
%
We may therefore apply Proposition~\ref{growth theorem} to $\partial_tu$ repeatedly on $(t_1,t_6)$ since
$t_6-t_1\ge 4T$ is assumed to be sufficiently large.
Therefore, there exist constants $h,\delta_3,c_3>0$ and integers $i_1,i_2$ with $1\le i_1\le i_2 \le j-1$ such that: if $1<i_1$, then we have
\begin{equation}
\label{growth1}
\begin{split}
&\text{either }\|\partial_tu\|_{L^2_{t,x}(t_1+iT,t_1+(i+1)T)} \le e^{-h T} \|\partial_tu\|_{L^2_{t,x} (t_1+(i-1)T,t_1+iT)}\\
&\text{or }{\delta_3}^{1/3}\|\partial_tu\|_{L^2_{t,x}(t_1+iT,t_1+(i+1)T)} \le
\|C_fe^{-2(t-t_0)}\|_{L^2_{t,x}(t_1+(i-1)T,\infty)}
\end{split}
\end{equation}
for every $i\in \{1,\dots,i_1-1\}$; if $i_1<i_2$, then we have
\begin{equation}
\label{growth2}
c_3 e^{-2(t-t_0)}\le \| \partial_tu(t) \|_{L^2_x}\le (3/2) \| \partial_tu(t') \|_{L^2_x}
\end{equation}
for every $t,t'\in(t_1+i_1T,t_1+i_2T)$ with $|t'-t|\le T$ and we have
\begin{equation}
\label{growth3}
\| \partial_t^2 u(t) \|_{L^2_x}\le (1/2) \| \partial_t u(t) \|_{L^2_x}
\end{equation}
for every $t\in (t_1+i_1T,t_1+i_2T)$; if $i_2<j-1$, then we have
\begin{equation}
\label{growth4}
\|\partial_t u\|_{L^2_{t,x} ( t_1+(i-1)T,t_1+iT) } \le e^{-h T}\|\partial_t u\|_{L^2_{t,x} (t_1+iT,t_1+(i+1)T)}
\end{equation}
for every $i\in \{i_2+1,\dots,j-1\}$.
Set $t_5=t_1+i_2T$.
Then, by \eqref{growth4} and \eqref{the first}, we have
\begin{equation}
\label{t_5 to t_6}
\begin{split}
\int_{t_5}^{t_6}\|\partial_t u(t)\|_{L^2_x}dt
&\le \sum_{i=i_2}^{j}\sqrt{T}\|\partial_t u\|_{L^2_{t,x} ( t_1+iT,t_1+(i+1)T)}\\
&\le \sqrt{T}(1-e^{-hT})^{-1}\sqrt{\delta}.
\end{split}
\end{equation}
In a similar way, by \eqref{growth1}, there exists a constant $C_{T,h}>0$ such that
\begin{equation}
\label{t_1 to t_2}
\int_{t_1}^{t_2}\|\partial_t u(t)\|_{L^2_x}dt\le C_{T,h}\sqrt{\delta}.
\end{equation}
If $i_1=i_2$, then by \eqref{t_0 to t_1} and \eqref{t_1 to t_2}, we have
\eqref{a-priori estimate};
we may therefore assume $i_1<i_2$.
Set $t_3=t_2+T/3$, $t_4=t_5-T/3$.
Then, by \eqref{the first}, we have
\begin{align}
\label{t_2 to t_3}
& \int_{t_2}^{t_3}\|\partial_t u(t)\|_{L^2_x}dt\le \sqrt{(T/3)\delta},
\; \int_{t_4}^{t_5}\|\partial_t u(t)\|_{L^2_x}dt\le \sqrt{(T/3)\delta}, \\
\label{t_2 to t_3'}
& \int_{t_2}^{t_3+T/4}\|\partial_t u(t)\|_{L^2_x}dt\le \sqrt{(7T/12)\delta}.
\end{align}
%

By \eqref{error2},
we may apply the Schauder estimate~\eqref{the Schauder estimate} to $w=\partial_tu,g=\partial_tf$.
Therefore, by \eqref{solution with 00} and \eqref{growth2}, there exists a constant $C_2>0$ such that for every $t\in[t_3,t_4]$, we have
\begin{equation}
\label{k>0}
\begin{split}
\|D_x\partial_t  u(t)\|_{L^2_x} 
\le C_2\|\partial_t u(t)\|_{L^2_x}.
\end{split}
\end{equation}
By \eqref{solution with 22'}, $u$ satisfies \eqref{Simon's equation} with $R$
satisfying \eqref{remainder estimate}.
Therefore, by \eqref{growth3} and \eqref{k>0},
for every $t\in[t_3,t_4]$, we have
\begin{equation*}
\| \partial_t u(t)+\grad E\bigl(u(t)\bigr) \|_{L^2_x}
= \| \partial_t^2u(t)+R\|_{L^2_x}
\le (3/4)\| \partial_t u(t) \|_{L^2_x}.
\end{equation*}
Therefore, by \eqref{solution with 22'} and \eqref{solution with 33},
we have \eqref{the assumption of Lojasiewicz-Simon}.
Therefore, by Proposition~\ref{Lojasiewicz-Simon}, we have
\begin{equation}
\label{u(t_3)0}
\int_{t_3}^{t_4}\|\partial_t u(t)\|_{L^2_x} dt \le (4/\theta) \left( \Bigl|E\bigl( u(t_3)\bigr)-E(0)\Bigr|^{\theta}+\delta^{\theta} \right)
\end{equation}
for some constant $\theta>0$.
Since $E$ satisfies \eqref{multiple integral} with \eqref{critical point}
and $u$ satisfies
the Schauder estimate~\eqref{the Schauder estimate},
there exist constants $C_3',C_3>0$ such that
\begin{equation}
\label{u(t_3)1}
\begin{split}
\bigl|E\bigl(u(t_3)\bigr)-E(0)\bigr|
&\le C_3'\|u(t_3)\|_{C^1_x}^2\\
&\le C_3\left(\sup_{t\in (t_3-T/4,t_3+T/4)}\|u(t)\|_{L^2_x}+e^{-2(t_3-t_0)}\right)^2.
\end{split}
\end{equation}
By \eqref{solution with 22}, \eqref{t_0 to t_1}, \eqref{t_1 to t_2} and \eqref{t_2 to t_3'},
there exists a constant $C_4>0$ such that
\begin{equation*}
\label{u(t_3)2}
\begin{split}
\sup_{t\in (t_3-T/4,t_3+T/4)}\|u(t)\|_{L^2_x}
\le \limsup_{t\to t_0}\|u(t)\|_{L^2_x}+ \int_{t_0}^{t_3+T/4}\|\partial_t u(t)\|_{L^2_x}dt
\le C_4\sqrt{\delta}.
\end{split}
\end{equation*}
By \eqref{growth2} and \eqref{the first}, there exists a constant $C_5>0$ such that
\begin{equation*}
\label{u(t_3)3}
e^{-2(t_3-t_0)}\le C_5\sqrt{\delta}.
\end{equation*}
Thus, \eqref{u(t_3)1} is bounded by $C_6\delta^{\theta}$ for some constant $C_6>0$.
Therefore, \eqref{u(t_3)0} is bounded by $C_7\delta^{\theta}$ for some constant $C_7>0$.
Therefore,
by \eqref{t_0 to t_1}, \eqref{t_5 to t_6}, \eqref{t_1 to t_2} and \eqref{t_2 to t_3},
we have \eqref{a-priori estimate}.
This completes the proof of Lemma~\ref{second lemma}.
\end{proof}
\section{Completion of the Proof}
\label{proof section}
We are ready now to complete the proof of Theorem \ref{key step}:
\begin{proof}[Proof of Theorem~$\ref{key step}$]
Let $\phi'$ be a parallel calibration of degree $m$ on the Euclidean space $(\mathbb{R}^n,g')$,
and let $\psi'$ be as in \eqref{psi'} in Section~\ref{statement section}, or equivalently as in \eqref{psi'2} in Section~\ref{lemma section 1}.
Let $X$ be a compact $\psi'$-submanifold of $S^{n-1}$.
Let $0<l<1$.
%
%
%
Suppose:
\begin{itemize}
\item[(S0)] $\epsilon>0$ is sufficiently small;
\item[(S1)] $0<a_0<b_0<a_1<b_1$, $a_0/b_0=a_1/b_1=l$;
\item[(S2)] $g$ is a Riemannian metric on $B(b_1)$ with
\[\|g-g'\|_{C^1}\le \epsilon,\|g-g'\|_{C^2}\le 1\]
with respect to $g'$; 
\item[(S3)] $\phi$ is a calibration of degree $m$ on $(B(b_1),g)$ with
\[(1+\log{\frac{b_1}{a_0}})\sup_{B(b_1)}|\phi-\phi'|\le \epsilon\]
where $|\bullet |$ is with respect to $g'$;
\item[(S4)] $M$ is a $\phi$-submanifold of $(\mathbb{R}^n,g)$, and $M$ is a closed subset of $(a_0,b_1)\times S^{n-1}$, where $(a_0,b_1)\times S^{n-1}$ is embedded into $\mathbb{R}^n$ by $(r,y)\mapsto ry$;
\item[(S5)] there exists a normal vector field $\nu_i$ on $(a_i,b_i)\times X$ in
$((a_i,b_i)\times S^{n-1},g')$, where $i=0,1,$ such that
\begin{equation*}
\begin{split}
M\cap ((a_i,b_i)\times S^{n-1})=G_{\mathrm{cyl}}(\nu_i)\text{ with }\| \nu _i \|_{C^1_{\mathrm{cyl}}}\le \epsilon
\end{split}
\end{equation*}
in the notation of Section~\ref{statement section}.
\end{itemize}
%
Let $\psi$ be as in \eqref{psi} in Section~\ref{lemma section 1}.
Then, by (S3) and (S5), we have
\begin{equation}
\label{psi estimate}
\sup_{r\in(a_0,b_0)\cup(a_1,b_1)}
\left|\int_{M\cap\{r\}\times S^{n-1}} \psi -\Vol(X)\right|\le C\epsilon
\end{equation}
for some constant $C>0$;
in the proof of Theorem~\ref{key step} a constant means a real number depending only on $l,m,n,X$ and $\phi'$.
%

By Proposition~\ref{dpsi-energy proposition}, the Stokes Theorem and \eqref{psi estimate}, we have
\begin{equation}
\label{energy estimate 2}
 \int_M|\mathrm{pr}_{TM^{\perp}}\partial_r|^2\dVol(M,g/r^2)\le  (2C/m)\epsilon+C_{m,n}\sup|\phi-\phi'|\Vol(M,g'/r^2).
 \end{equation}
Therefore, by Proposition~\ref{vol estimate} and \eqref{psi estimate}, we have
\begin{equation*}
\begin{split} \Vol(M,g'/r^2)
&\le C'\log{\frac{b_1}{a_0}}
+C'(1+m\log{\frac{b_1}{a_0}}) \left(\epsilon+\sup|\phi-\phi'|\Vol(M,g'/r^2)\right)
\end{split}
\end{equation*}
for some constant $C'>0$. Therefore, by (S3), we have
%
\begin{equation}
\label{vol estimate'}
 \Vol(M,g'/{r^2})\le C''\log{\frac{b_1}{a_0}}
\end{equation}
for some constant $C''>0$.
Therefore, by \eqref{energy estimate 2}, we have
\begin{equation}
\label{energy small}
\int_M|\mathrm{pr}_{TM^{\perp}}\partial_r|^2\dVol(M,g/r^2)\le C'''\epsilon
\end{equation}
for some constant $C'''>0$.
%
%
%
%
Choose a constant $\epsilon_*>0$ so that if $I$ is an open interval of $(0,\infty)$,
and if $\nu$ is a normal vector field
on $I\times X$ in $(I\times S^{n-1},g')$ with $\|\nu\|_{C^0_{\mathrm{cyl}}}\le \epsilon_*$, then
$G_{\mathrm{cyl}}(\nu)$ is contained in a tubular neighbourhood of $I\times X$ in $(I\times S^{n-1},g')$.
Here, $G_{\mathrm{cyl}}(\nu)$ is as in Section~\ref{statement section}.
If $\epsilon_*$ is sufficiently small, then we have
\[|r\partial_r(\nu/r)|^2 \le 2 |\mathrm{pr}_{TM^{\perp}}\partial_r|^2,\]
as in \cite[(7.13), p561]{Simon} or \cite[3.2, Part~I]{Simon2}.
Therefore, by \eqref{energy small}, we have
\begin{equation}
\label{L^2-estimate}
\int_{M\cap(I\times S^{n-1})}|r\partial_r(\nu/r)|^2\dVol(M,g/r^2) \le 2 C'''\epsilon.
\end{equation}
%
%
Choose $0<\lambda<\lambda''<\lambda'<1$ so that $l\lambda'<\lambda<\lambda''<l<\lambda'$.
By \eqref{energy small}, we may apply Lemma~\ref{first lemma} to
$M\cap((\lambda b_1,b_1))\times S^{n-1}$.
Therefore, there exists a normal vector field $\nu$ on $(\lambda b_1,b_1)\times X$ in
$((\lambda b_1,b_1)\times S^{n-1},g'/r^2)$
such that
\begin{equation}
\label{non-empty}
\begin{split}
& M\cap ((\lambda b_1,b_1)\times S^{n-1})=G_{\mathrm{cyl}}(\nu)
\text{ with } \|\nu|_{(\lambda'' b_1,\lambda' b_1)\times X}\|_{C^{1,1/2}_{\mathrm{cyl}}} \le \epsilon_*.
\end{split}
\end{equation}
Let $S_*$ be the set of all $b_*\in [\lambda' a_0/\lambda,a_1)$ such that
there exists a normal vector field $\nu$ on $(b_*,a_1)\times X$ in $((b_*,a_1)\times S^{n-1},g'/r^2)$
such that
\begin{equation}
\label{nu}
\begin{split}
& M\cap ((b_*,b_1)\times S^{n-1})=G_{\mathrm{cyl}}(\nu)
\text{ with } \|\nu|_{(b_*,a_1)\times X}\|_{C^{1,1/2}_{\mathrm{cyl}}} \le \epsilon_*.
\end{split}
\end{equation}
$S_*$ is non-empty since $\lambda''b_1\in S_*$ by \eqref{non-empty}.
%
%
\begin{prop}
\label{S_*}
Suppose $b_*\in S_*\cap [\lambda a_0/\lambda',b_1)$,
and let $\nu$ be as in \eqref{nu}.
Then, there exist constants $c_{10},C_{10}>0$
such that
\[ \|\nu|_{(b_*,b_1)\times X}\|_{C^1_{\mathrm{cyl}}}
\le C_{10}\epsilon^{c_{10}}.\]
\end{prop}
\begin{proof}
By Proposition~\ref{psi minimal}, $X$ is a minimal submanifold of $S^{n-1}$.
By (S4), $M$ is a minimal submanifold of $((a_0,b_1)\times S^{n-1},g)$ with $\|g-g'\|_{C^2(B(b_1))}\le 1$ as in (S2).
Set
\[u(t,x)=e^{t/m}\nu(e^{-t/m}x),\;t_0=-m\log{b_1},t_*=-m\log{b_*}.\]
Then, by a result of Simon~\cite[Remark~3.3, Part~I]{Simon2}, $u$ satisfies \eqref{Simon's equation}
for some $E,R,f$ depending only on $m,n,X.$
We shall apply Lemma~\ref{second lemma} to $u$.
By \eqref{nu}, we have
\begin{equation}
\label{giving solution with 22}
\|\nu|_{(b_*,a_1)\times X}\|_{C^{1,1/2}_{\mathrm{cyl}}}
=\| u \|_{C^{1,1/2}_{t,x}(-m\log{a_1},t_*)}\le \epsilon_*.
\end{equation}
Therefore, we have \eqref{solution with 11}.
By (S5), we have \eqref{solution with 22}.
By \eqref{L^2-estimate}, \eqref{giving solution with 22} and (S2),
there exists a constant $C_{11}>0$ such that
\begin{equation*}
\begin{split}
&\|\partial_tu\|_{L^2_{t,x}(t_0,t_*)}^2
=\int_{(b_*,b_1)\times X} |r\partial_r(\nu/r)|^2 dr/r\dVol(X)
\le C_{11}{\epsilon}.
\end{split}
\end{equation*}
Therefore, we have \eqref{solution with 44}.
It suffices therefore to prove \eqref{solution with 33}.
In a way similar to \eqref{psi estimate}, by \eqref{nu}, we have
\begin{equation*}
\label{giving solution with 33}
\begin{split}
&\sup_{b\in (b_*,b_1)} \Vol(X) -\Vol \left( M\cap \{b\}\times S^{n-1},g'/{r^2} \right)\\
&\le \sup_{b\in (b_*,b_1),b'\in(a_1,b_1)} \int_{ M\cap \{ b'\}\times S^{n-1}}\psi -\int_{M\cap \{b\}\times S^{n-1}}\psi+C_{12}\epsilon
\end{split}
\end{equation*}
for some constant $C_{12}>0$.
By Proposition~\ref{dpsi-energy proposition}, \eqref{energy small}, \eqref{vol estimate'} and (S3), we have
\begin{equation*}
\begin{split}
\sup_{b\in (b_*,b_1),b'\in(a_1,b_1)} \int_{M\cap \{ b'\}\times S^{n-1} }\psi -\int_{M\cap \{ b\}\times S^{n-1}}\psi
\le C_{13}{\epsilon}
\end{split}
\end{equation*}
for some constant $C_{13}>0$.
Thus, there exists a constant $C_{14}>0$ such that
\begin{equation*}
\sup_{b\in (b_*,b_1)} \Vol(X) -\Vol \left( M\cap \{b\}\times S^{n-1},g'/{r^2} \right)\le C_{14}{\epsilon}.
\end{equation*}
Therefore, we have \eqref{solution with 33}.
We may now apply Lemma~\ref{second lemma} to $u$.
Therefore, as in \eqref{conclusion of second lemma}, we have
\[\sup_{t\in(t_0,t_*)}\|u\|_{L^{2}_{x}}\le C_{15}\epsilon^{c_{15}}\]
for some constants $c_{15},C_{15}>0$.
Therefore, by interpolation and \eqref{giving solution with 22}, we have
\begin{equation*}
\begin{split}
\|\nu|_{(b_*,a_1)\times X}\|_{C^1_{\mathrm{cyl}}}
=\| u \|_{C^1_{t,x}(-m\log{a_1},t_*)}
\le C_{10}\epsilon^{c_{10}}
\end{split}
\end{equation*}
for some constants $C_{10},c_{10}>0$.
By (S5), this proves Proposition~\ref{S_*}.
\end{proof}
%
Suppose $b_*\in S_*$.
Then, by Proposition~\ref{S_*}, we may apply Lemma~\ref{first lemma} to
$M\cap((\lambda b_*/\lambda',b_*/\lambda')\times S^{n-1})$. 
Therefore,
$\lambda''b_*/\lambda'\in S_*$.
Therefore, $b_*$ is an interior point in $S_*$.
$S_*$ is thus an open subset of $[\lambda' a_0/\lambda,a_1)$.

By definition, $S_*$ is a closed subset of $[\lambda' a_0/\lambda,a_1)$.
$S_*$ is thus a non-empty open closed subset of $[\lambda' a_0/\lambda,a_1)$.
Therefore, $S_*=[\lambda' a_0/\lambda,a_1)$.
Therefore, $\lambda' a_0/\lambda\in S_*$.
Therefore, by \eqref{nu} and Proposition~\ref{S_*}, we have
\begin{equation*}
M\cap ((\lambda' a_0/\lambda,b_1)\times S^{n-1})=G_{\mathrm{cyl}}(\nu)
\text{ with } \|\nu|_{(\lambda' a_0/\lambda,b_1)\times X}\|_{C^{1}_{\mathrm{cyl}}} \le C_{10}\epsilon^{c_{10}}.
\end{equation*}
Therefore, by (S5), we have \eqref{conclusion of key step}. This completes the proof of Theorem~\ref{key step}.
\end{proof}

\end{document}